\documentclass[a4paper,12pt]{article}

\usepackage{amsmath}
\usepackage{amssymb}
\usepackage{amsthm}

\usepackage[english]{babel}

\usepackage[T2A]{fontenc}
\usepackage{graphics}
\usepackage{epsfig}

\newtheorem{theorem}{Theorem}

\newtheorem{lemma}{Lemma}

\theoremstyle{remark}

\theoremstyle{definition}

\theoremstyle{definition}

      \newcommand {\al}   {\alpha}          
                
      \newcommand {\del}  {\delta}          
              \newcommand {\ve}   {\varepsilon}
                 
      \newcommand {\lam}  {\lambda}         
      
                \newcommand {\Om}  {\Omega}
      \newcommand {\pl}   {\partial}        
           
      \newcommand {\RRR}  {{\mathbb R}}     
           \newcommand {\GGG}  {{\cal G}}
              \newcommand {\BBB}  {{\cal B}}
              \newcommand {\FFF}  {{\cal F}}

              \newcommand {\BBBB}  {B}
     \newcommand {\beq}  {\begin{equation}}
      \newcommand {\eeq}  {\end{equation}}  \newcommand {\lab}  {\label}

\author{Alexander Plakhov\thanks{Aberystwyth University, Aberystwyth SY23 3BZ, UK, on leave from Department of Mathematics, University of Aveiro, Aveiro 3810-193, Portugal}
\and Alena Aleksenko\thanks{Department of Mathematics, Aveiro
University, Aveiro 3810, Portugal}
 \thanks{This work was supported by {\it Centre for Research on Optimization
and Control} (CEOC) from the ''{\it Funda\c{c}\~{a}o para a
Ci\^{e}ncia e a Tecnologia}'' (FCT), cofinanced by the European
Community Fund FEDER/POCTI, and by the FCT research project
PTDC/MAT/72840/2006. }}

\title{The problem of the body of revolution of minimal resistance}
\date{}

\begin{document}

\maketitle

\begin{abstract}
Newton's problem of the body of minimal aerodynamic resistance is traditionally
stated in the class of {\it convex} axially symmetric bodies with
fixed length and width. We state and solve the minimal resistance
problem in the wider class of axially symmetric but {\it generally
nonconvex} bodies. The infimum in this problem is not attained. We
construct a sequence of bodies minimizing the resistance. This
sequence approximates a convex body with smooth front surface, while
the surface of approximating bodies becomes more and more
complicated. The shape of the resulting convex body and the value of
minimal resistance are compared with the corresponding results for
Newton's problem and for the problem in the intermediate class of
axisymmetric bodies satisfying the {\it single impact} assumption
\cite{CL1}. In particular, the minimal resistance in our class is
smaller than in Newton's problem; the ratio goes to $1/2$ as
(length)/(width of the body) $\to 0$, and to $1/4$ as
(length)/(width) $\to +\infty$.
\end{abstract}

\begin{quote}
\end{quote}

\begin{quote}
\end{quote}

\begin{quote}
{\small {\bf Running title:} Problem of minimal resistance}
\end{quote}

\section{Introduction}

In 1687, I. Newton in his {\it Principia} \cite{Principia}
considered a problem of minimal resistance for a body moving in a
homogeneous rarefied medium. In slightly modified terms, the problem
can be expressed as follows.

A convex body is placed in a parallel flow of point particles. The
density of the flow is constant, and velocities of all particles are
identical. Each particle incident on the body makes an elastic
reflection from its boundary and then moves freely again. The flow
is very rare, so that the particles do not interact with each other.
Each incident particle transmits some momentum to the body; thus,
there is created a force of pressure on the body; it is called {\it
aerodynamic resistance force}, or just {\it resistance}.

\begin{figure}\label{newtso}
\begin{center}
\includegraphics[width=0.30\textwidth, angle=0
]{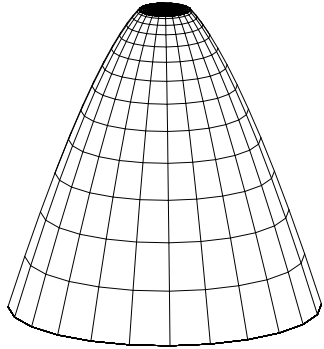} \caption{The Newton solution for $h =
2$.}\label{newtsol2}
\end{center}
\end{figure}

Newton described (without proof) the body of minimal resistance in
the class of convex and axially symmetric bodies of fixed length and
maximal width, where the symmetry axis is parallel to the flow
velocity. That is, any body from the class is inscribed in a right
circular cylinder with fixed height and radius. The rigorous proof
of the fact that the body described by Newton is indeed the
minimizer was given two centuries later. From now on, we suppose
that the radius of the cylinder equals 1 and the height equals $h$,
with $h$ being a fixed positive number. The cylinder axis is
vertical, and the flow falls vertically downwards. The body of least
resistance for $h = 2$ is shown on fig.\,\ref{newtso}.

Since the early 1990s, there have been obtained new interesting
results related to the problem of minimal resistance in various
classes of admissible bodies \cite{BK}-\cite{P2}. In particular,
there has been considered the wider class of convex (generally
non-symmetric) bodies inscribed in a given cylinder
\cite{BK}-\cite{BrFK},\cite{LP2},\cite{LO}. It was shown that the
solution in this class exists and does not coincide with the Newton
one. The problem is not completely solved till now. The numerical
solution for $h = 1.5$ is shown on fig.\,2.\footnote{This figure is
reproduced with kind permission of E. Oudet.}

\begin{figure}\label{nonnewtso}
\begin{center}
\includegraphics[height=5cm, width=5cm,angle=90]{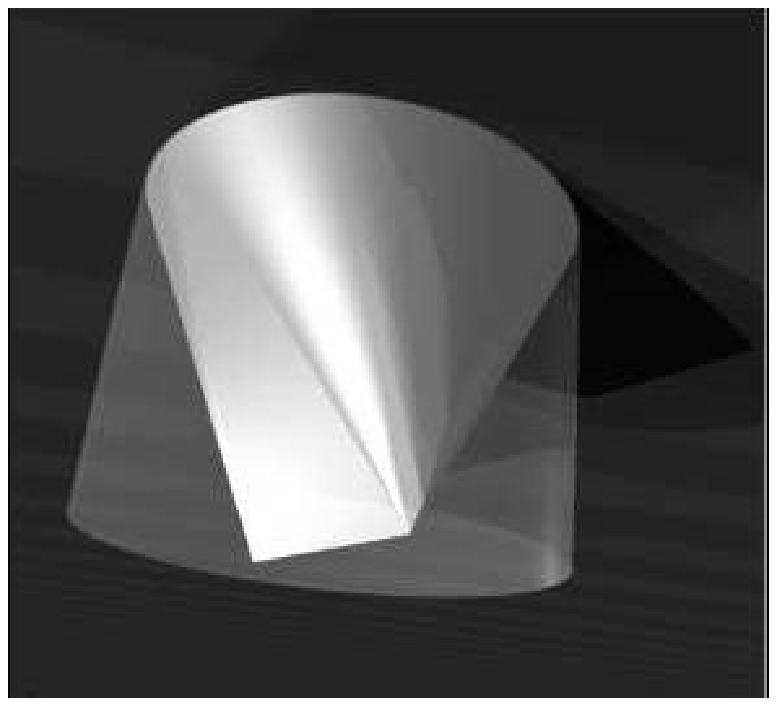}
\caption{The non-symmetric solution for $h =1.5$.}\label{newtsol}
\end{center}
\end{figure}


By removing both assumptions of symmetry and convexity, one gets the
(even wider) class of bodies inscribed in a given cylinder. More
precisely, a generic body from the class is a connected set with
piecewise smooth boundary which is contained in the cylinder,
contains an orthogonal cross section of the cylinder, and satisfies
a regularity condition to be specified below. Notice that there may
occur multiple reflections of particles from the surface of a {\it
non-convex} body, while reflections from {\it convex} bodies are
always single. The problem of minimal resistance in this class was
solved in \cite{P1,P2}. In contrast to the class of convex $\&$
axisymmetric bodies and the class of convex bodies, the infimum of
resistance here equals zero, and we believe the infimum cannot be
attained.

In addition to the classes of admissible bodies discussed above:\\
\hspace*{6.5mm}(i) convex $\&$ axisymmetric (the classical Newton problem);\\
\hspace*{6.5mm}(ii) convex but generally non-symmetric;\\
\hspace*{6.5mm}(iii) generally nonconvex and non-symmetric,\\
there remains a class that has not been studied as yet:\\
\hspace*{6.5mm}(iv) axisymmetric but generally nonconvex bodies.\\
The aim of this paper is to fill this gap: we shall solve the
minimal resistance problem for the fourth class.

Note that in the paper \cite{CL1} there was considered the intermediate class of\\
\hspace*{6.5mm}(v) axially symmetric nonconvex bodies, under the additional so-called\\
\hspace*{11mm} "single impact assumption".\\
This geometric assumption on the body's shape means that each
particle hits the body at most once; multiple reflections
are not allowed. On the contrary, multiple reflections are allowed in
our setting; we only assume that the body's boundary is piecewise
smooth and satisfies the regularity condition stated below.

The class (v) is intermediate between the classes (i) and (iv); it
contains the former one and is contained in the latter one. We shall
determine the minimal resistance and the minimizing sequence of
bodies for the class (iv) (which will be referred to as {\it
nonconvex case}), and compare them with the corresponding results
for the class (i) ({\it Newton} case) and for the class (v) ({\it
single impact} case).\footnote{Note that Newton himself did not
state explicitly the assumption of convexity; in this sense, the
cases (iv) and (v) can be regarded as "relaxed versions"' of the
Newton problem.}

Consider a compact connected set $\BBBB \subset \RRR^3$ and choose
an orthogonal reference system $Oxyz$ in such a way that the axis $Oz$ is
parallel to the flow direction; that is, the particles move
vertically downwards with the velocity $(0, 0, -1)$. Suppose that a
flow particle (or, equivalently, a billiard particle in $\RRR^3
\setminus \BBBB$) with coordinates $x(t) = x$,\, $y(t) = y$,\, $z(t)
= -t$ makes a finite number of reflections at regular points of the
boundary $\pl \BBBB$ and moves freely afterwards. Denote by
$\nu_\BBBB(x,y)$ the final velocity. If there are no reflections,
put $\nu_\BBBB(x,y) = (0, 0, -1)$.

Thus, one gets the function $\nu_\BBBB = (\nu_\BBBB^x, \nu_\BBBB^y,
\nu_\BBBB^z)$ taking values in $S^2$ and defined on a subset of
$\RRR^2$. We impose the regularity condition requiring that $\nu_\BBBB$ {\it is
defined on a full measure subset of} $\RRR^2$. All convex sets
$\BBBB$ satisfy this condition; examples of non-convex sets violating
it are given on figure 3. Both sets are of the form $B =
G \times [0,\, 1] \subset \RRR^2_{x,z} \times R^1_y$, with $G$ being
shown on the figure. On fig.\,3a, a part of the boundary is an arc
of parabola with the focus $F$ and with the vertical axis. Incident particles, after making a reflection from the arc, get into
the singular point $F$ of the boundary. On fig.\,3b, one part of the
boundary belongs to an ellipse with foci $F_1$ and $F_2$, and
another part, $AB$, belongs to a parabola with the focus $F_1$ and with the
vertical axis. After reflecting from $AB$, particles of the flow get
trapped in the ellipse, making infinite number of reflections and
approaching the line $F_1 F_2$ as time goes to $+\infty$. In both
cases, $\nu_\BBB$ is not defined on the corresponding
positive-measure subsets of $\RRR^2$.

\begin{figure}\label{sing}
\begin{center}
\includegraphics[height=10cm,angle=-90]{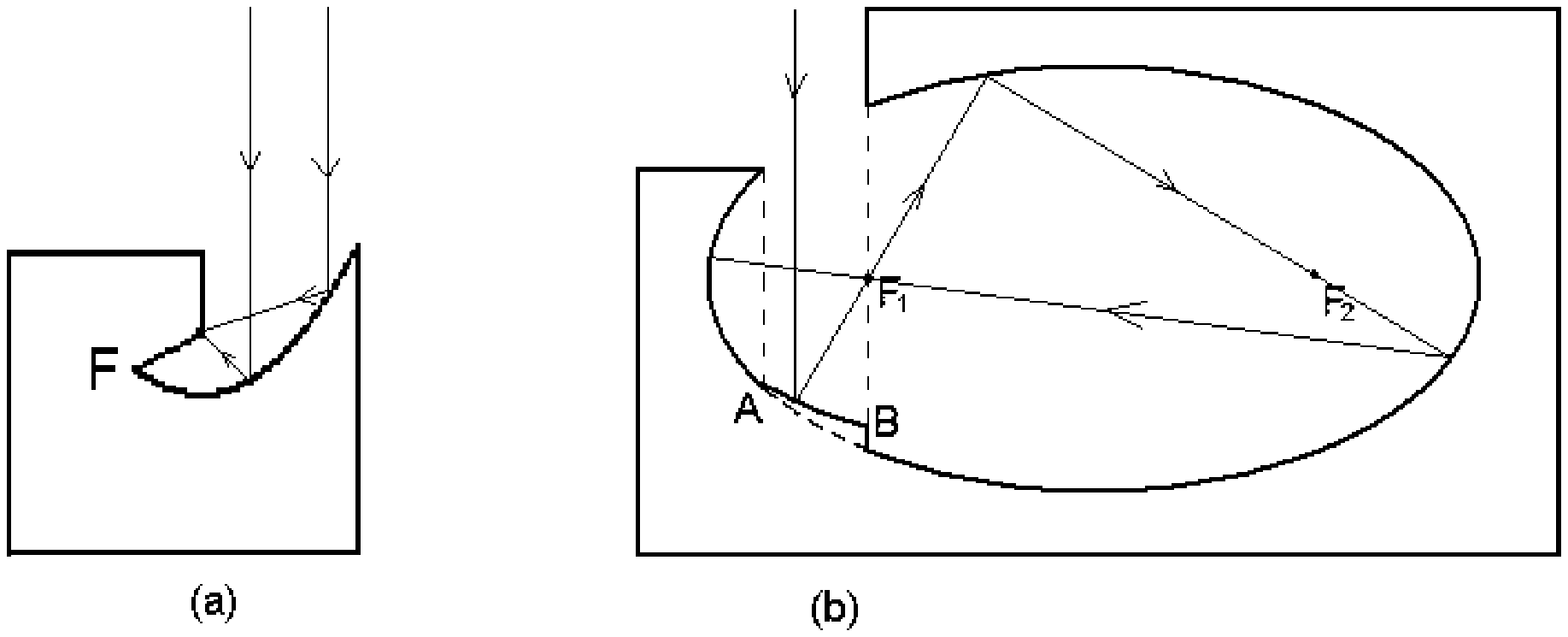}
\caption{(a) After reflecting from the arc of parabola, the
particles get into the singular point $F$. (b) After reflecting from
the arc of parabola $AB$, the particles get trapped in the
ellipse.}\label{sing2}
\end{center}
\end{figure}



Each particle interacting with the body $\BBBB$ transmits to it the momentum
equal to the particle mass times $\left( (0, 0, -1) -
\nu_\BBBB(x,y) \right)$. Summing up over all momenta transmitted per
unit time, one obtains that the resistance of $\BBBB$
equals $-\rho\, \mathrm{R}(\BBBB)$, where
\begin{equation*}\label{resist general formula}
\mathrm{R}(\BBBB) = \int\!\!\!\!\int_{\RRR^2} \left( \nu_\BBBB^x,\,
\nu_\BBBB^y,\, 1 + \nu_\BBBB^z \right) dx\, dy,
\end{equation*}
and $\rho$ is the flow density. One is usually
interested in minimizing the third component of
$\mathrm{R}(\BBBB)$,\footnote{Note that in the axisymmetric cases
(i), (iv), and (v), the first and second components of $\mathrm{R}(\BBBB)$
are zeros, due to radial symmetry of the functions $\nu_\BBBB^x$ and
$\nu_\BBBB^y:\, \ \mathrm{R}_x(\BBBB) = 0 = \mathrm{R}_y(\BBBB)$.}
\begin{equation}\label{R_z general formula}
\mathrm{R}_z(\BBBB) = \int\!\!\!\!\int_{\RRR^2} \left( 1 +
\nu_\BBBB^z(x,y) \right) dx\, dy.
\end{equation}

If $\BBBB$ is convex then the upper part of the boundary $\pl \BBBB$
is the graph of a concave function $w(x,y)$. Besides, there is at
most one reflection from the boundary, and the velocity of the
reflected particle equals $\nu_\BBBB(x,y) = (1 + |\nabla w|^2)^{-1}
(-2w_x,\, -2w_y,\, 1 - |\nabla w|^2)$. Therefore, the formula
(\ref{R_z general formula}) takes the form
\begin{equation}\label{R_z for convex body}
\mathrm{R}_z(\BBBB) = \int\!\!\!\!\int \frac{2}{1 + |\nabla
w(x,y)|^2}\ dx\, dy,
\end{equation}
the integral being taken over the domain of $w$.

Further, if $\BBBB$ is a convex axially symmetric body then (in a
suitable reference system) the function $w$ is radial: $w(x,y) =
f(\sqrt{x^2 + y^2})$, therefore one has
\begin{equation}\label{R_z convex radial}
\mathrm{R}_z(\BBBB) = 2\pi\! \int\! \frac{2r}{1 + f'^2(r)}\ dr,
\end{equation}
the integral being taken over the domain of $f$.

Thus, in the cases (i), (ii), and (v) the problem of minimal
resistance reads as follows:
 \beq\label{sluchaj (i)}
\text{(i) \, \ \ \ \ \ \ \  minimize} \ \ \ \ \int_0^1 \frac{r}{1 +
f'^2(r)}\ dr \hspace*{55mm}
 \eeq
 \vspace{-5mm}

over all concave monotone non-increasing functions $f: [0,\, 1] \to
[0,\, h]$;
$$
\text{(ii) \ \ \ minimize} \ \ \int\!\!\!\!\int_\Om\, \frac{1}{1 +
|\nabla w(x,y)|^2}\ dx\, dy \hspace*{50mm}
$$
 \vspace{-5mm}

over all concave functions $w: \Om \to [0,\, h]$, where $\Om = \{
x^2 + y^2 \le 1 \}$ is the

unit circle;
$$
\text{(v) minimize the functional (\ref{sluchaj (i)}) over the set
$\mathcal{C}_h$ of functions $f : [0,\, 1] \to [0,\, h]$}
$$

 \vspace{-3mm}

satisfying the single impact condition (see \cite{CL1}, formulas (3)
and (1)).
 \vspace{2mm}

In the nonconvex cases (iii) and (iv) the functional to be minimized
(\ref{R_z general formula}) cannot be written down explicitly in
terms of the body's shape. Still, in the radial case (iv) it can be
simplified in the following way.

Let $\BBBB$ be a compact connected set inscribed in the cylinder
$x^2 + y^2 \le 1$,\, $0 \le z \le h$ and possessing rotational
symmetry with respect to the axis $Oz$. This set is uniquely defined
by its vertical central cross section $G = \{ (x,z): (x,0,z) \in
\BBBB \}$. It is convenient to reformulate the problem in terms of
the set $G$.

\begin{figure}\label{scat}
\begin{center}
\includegraphics[width=5cm, angle=0]{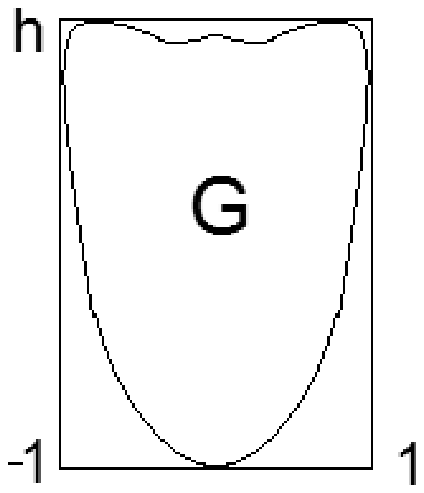}
\caption{A set $G \in \GGG_h$.}\label{newt1}
\end{center}
\end{figure}

Consider the billiard in $\RRR^2 \setminus G$ and suppose that a
billiard particle initially moves according to $x(t) = x$,\, $z(t) =
-t$, then makes a finite number of reflections (maybe none) at
regular points of $\pl G$, and finally moves freely with the
velocity $v_G(x) = (v_G^x(x), v_G^z(x))$. The regularity condition
now means that that the so determined function $v_G$ is defined for
almost every $x$. One can see that $\nu_\BBBB^x(x,y) =
(x/\!\sqrt{x^2 + y^2}) v_G^x(\sqrt{x^2 + y^2})$,\, $\nu_\BBBB^y(x,y)
= (y/\!\sqrt{x^2 + y^2}) v_G^y(\sqrt{x^2 + y^2})$, and
$\nu_\BBBB^z(x,y) = v_G^z(\sqrt{x^2 + y^2})$. It follows that
$\mathrm{R}_x(\BBBB) = 0 = \mathrm{R}_y(\BBBB)$ and
$\mathrm{R}_z(\BBBB) = 2\pi \int_0^1 (1 + v_G^z(x))\, x\, dx$. Thus,
our minimization problem takes the form
\begin{equation}\label{the problem}
\inf_{G \in \GGG_h} R(G), \ \ \ \text{where} \ \ R(G) = \int_0^1\!
\left( 1 + v_G^z(x) \right) x\, dx
\end{equation}
and $\GGG_h$ is the class of compact connected sets $G \subset
\RRR^2$ with piecewise smooth boundary that are inscribed in the
rectangle $-1 \le x \le 1$,\, $0 \le z \le h$,\footnote{That is,
belong to the rectangle and have nonempty intersection with each
of its sides.} are symmetric with respect to the axis $Oz$, and
satisfy the regularity condition (see fig. \ref{newt1}).

The main results are stated in section 2: the minimization problem
is solved and the solution is compared with the Newton solution
(case (i)) and the single-impact solution (case (v)). Details of all
proofs are put in section 3.

\section{Statement of the results}

\begin{figure}\label{scat2}
\begin{center}
\includegraphics[width=0.30\textwidth, angle=0]{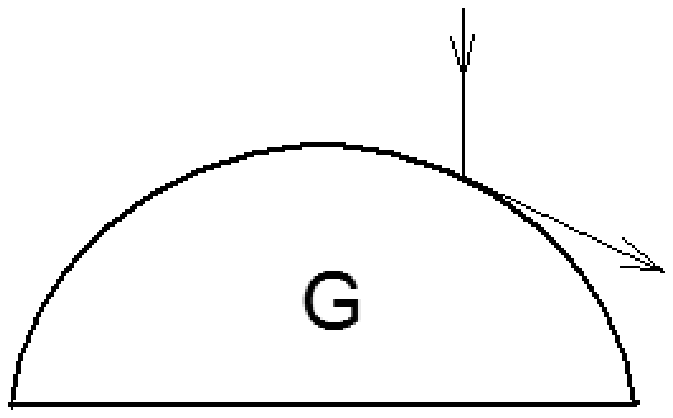}
\caption{Modified reflection law.}\label{newt2}
\end{center}
\end{figure}

Denote by $\GGG_h^{\text conv}$ the class of convex sets from
$\GGG_h$. One can easily see that if $G \in \GGG_h$ then  conv$\,G
\in \GGG_h^{\text conv}$. For $G \subset \GGG_h^{\text conv}$ define
the {\it modified law} of reflection as follows. A particle
initially moves vertically downwards according to $x(t) = x$,\,
$z(t) = -t$ and reflects at a regular point of the boundary $\pl G$;
at this point the velocity instantaneously changes to $\hat{v}_G(x)
= (\hat{v}_G^x(x), \hat{v}_G^z(x))$, where $\hat{v}_G(x)$ is the
unit vector tangent to $\pl G$ such that $\hat{v}_G^z(x) \le 0$
and $x \cdot \hat{v}_G^x(x) \ge 0$ (see fig.\,\ref{newt2}).

The set $G \in \GGG_h$ is bounded above by the graph of a concave
even function $z = f_G(x)$. For $x > 0$, one has
 \beq\label{vel}
\hat{v}_G(x) = \frac{(1,f'_G(x))}{\sqrt{1 +
f'^2_G(x)}}\,.
 \eeq
The resistance of $G$ under the modified reflection law equals $(0,
-\widehat R(G))$, where
\begin{equation}\label{respl}
\widehat{R}(G) = \int_0^1 (1 + \hat{v}^z_G(x))\, x\, dx\,.
\end{equation}
Taking into account (\ref{vel}), one gets
 \beq\label{resmod}
\widehat R(G) = \int_0^1 \left( 1 + \frac{f'_G(x)}{\sqrt{1 +
f'^2_G(x)}} \right) x\, dx;
 \eeq
the function $f_G$ is concave, nonnegative, and monotone
non-increasing, with  $f(0) = h$.

\begin{theorem}\label{t1}
\begin{equation}\label{inf inf}
\inf_{G \in \GGG_h} R(G) = \inf_{G \in \GGG_h^{\text conv}}
\widehat{R}(G).
\end{equation}
\end{theorem}

This theorem follows from the following lemmas \ref{l1} and \ref{l2}
which will be proved in the next section.

\begin{lemma}\label{l1}
For any $G \in \mathcal G_h$ one has
\begin{equation*}
R(G) \ge \widehat{R}(\text{conv}\, G).
\end{equation*}
\end{lemma}

\begin{lemma}\label{l2}
Let $G \in \GGG_h^{\text conv}$. Then there exists a sequence of
sets $G_n \in \GGG_h$ such that
$$
\lim_{n\to\infty} R(G_n) = \widehat R(G).
$$
\end{lemma}

Indeed, lemma \ref{l1} implies that $\inf_{G \in \GGG_h} R(G) \ge
\inf_{G \in \GGG_h^{\text conv}} \widehat{R}(G)$, and lemma \ref{l2}
implies that $\inf_{G \in \GGG_h} R(G) \le \inf_{G \in \GGG_h^{\text
conv}} \widehat{R}(G)$.

Theorem \ref{t1} allows one to state the minimization problem
(\ref{the problem}) in an explicit form. Namely, taking into account
(\ref{resmod}) and putting $f = h - f_G$, one rewrites the right
hand side of (\ref{inf inf}) as
\begin{equation}\label{probl}
\inf_{f \in \FFF_h} \int_0^1 \left( 1 - \frac{f'(x)}{\sqrt{1 +
f'^2(x)}} \right) x\, dx,
\end{equation}
where $\FFF_h$ is the set of convex monotone non-decreasing
functions $f: [0,\, 1] \to [0,\, h]$ such that  $f(0) = 0$. The
solution of (\ref{probl}) is provided by the following general
theorem.

Consider a positive piecewise continuous function $p$ defined on
$\RRR_+ := [0,\, +\infty)$ and converging to zero as $u \to
+\infty$, and consider the problem
 \beq\label{1}
\inf_{f \in \FFF_h} \mathcal{R}[f],~~~ \text{where}~~~
\mathcal{R}[f] = \int_0^1 p(f'(x))\, x\, dx.
 \eeq
Denote by  $\bar p(u)$,\, $u \in \RRR_+$ the greatest convex
function that does not exceed $p(u)$. Put $\xi_0 = -1/\bar p'(0)$
and $u_0 = \inf \{ u > 0:\, \bar p(u) = p(u) \}$. One always has
$\xi_0 \ge 0$; if $u_0 = 0$ and there exists $p'(0)$ then $\xi_0 =
-1/p'(0)$, and if  $u_0 > 0$ then $\xi_0 = u_0/(p(0) - p(u_0))$.
Denote by $u = \upsilon(z)$,\, $z \ge \xi_0$ the generalized inverse
of the function $z = -1/\bar p'(u)$, that is, $\upsilon(z) = \sup \{
u : -1/\bar p'(u) \le z \}$. By $\Upsilon$, denote the primitive of
$\upsilon$:\, $\Upsilon(z) = \int_{\xi_0}^z \upsilon(\xi) d\xi$,\,
$z \ge \xi_0$. Finally, put $\mathcal{R}(h) := \inf_{f \in \FFF_h}
\mathcal{R}[f]$.

\begin{theorem}\label{t2}
For any $h > 0$ the solution $f_h$ of the problem (\ref{1}) exists
and is uniquely determined by
 \beq\label{f solution}
f_h(x) = \left\{ \begin{array}{ll} 0 & \text{if}~~ 0 \le x \le x_0\\
\frac{1}{z}\, \Upsilon(z x) & \text{if}~~ x_0 < x \le 1\,,
 \end{array} \right.
 \eeq
where $z = z(h)$ is a unique solution of the equation
 \beq\label{z solution}
\Upsilon(z) = z h
 \eeq
and $x_0 = x_0(h) = \xi_0/z(h)$. Further, one has $f_h'(x_0 + 0) =
u_0$. The function $x_0(h)$ is continuous and $x_0(0) = 1$. The
minimal resistance equals
 \beq\label{7}
\mathcal{R}(h)\ =\ \frac 12 \left ( \bar p(\upsilon(z))\, +\,
\frac{\upsilon(z) - h}{z} \right );
 \eeq
in particular, $\mathcal{R}(0) = p(0)/2$.

If, additionally, the function $p$ satisfies the asymptotic relation
$p(u) = cu^{-\al}\, (1 + o(1))$ as $u \to +\infty$,\, $c > 0$,\,
$\al > 0$ then
 \beq\label{9}
x_0(h) = c\al \left( \frac{\al + 1}{\al + 2} \right)^{\al + 1} \xi_0
h^{-\al - 1} (1 + o(1)),~~  h \to +\infty,
 \eeq
 and
 \beq\label{8}
\mathcal{R}(h) = \frac{c}{2} \left( \frac{\al + 1}{\al + 2}
\right)^{\al + 1} h^{-\al} (1 + o(1)),~~ h \to +\infty.
 \eeq
\end{theorem}

Let us apply the theorem to the three cases under consideration.
 \vspace{1mm}

{\bf 1.} First consider the {\it non-convex case}. The problem
(\ref{probl}) we are interested in is a particular case of (\ref{1})
with $p(u) = p_{\text nc}(u) := 1 - u/\sqrt{1 + u^2}$ (the subscript
"nc" stands for "non-convex"). The function $p_{\text nc}$ itself, however, is
convex, hence $u_0 = 0$ and $\bar p_{\text nc} \equiv p_{\text nc}$.
Further, one has $-1/\bar p'_{\text nc}(u) = (1 + u^2)^{3/2}$,
therefore $\upsilon_{\text nc}(z) = \sqrt{z^{2/3} - 1}$,\,
$\xi_0^{\text nc} = 1$, and
 \beq\label{gamma}
\Upsilon_{\text nc}(z) = \frac 38 (2z^{2/3} - 1) z^{1/3}
\sqrt{z^{2/3} - 1} - \frac{3}{8} \ln(z^{1/3} + \sqrt{z^{2/3} - 1}).
 \eeq
The formulas  (\ref{gamma}), (\ref{z solution}), and (\ref{f
solution}) with $x_0 = 1/z$, determine the solution of
(\ref{probl}). Notice that, as opposed to the Newton case, the
solution is given by the {\it explicit} formulas. However, they
contain the parameter $z$ to be defined implicitly from (\ref{z
solution}).

Further, according to theorem 2, $f'_h(x_0 + 0) = 0 = f'_h(x_0 -
0)$,\, $x_0 = x_0^{\text nc}$, hence the solution $f_h$ is
differentiable everywhere in $(0,\, 1)$. Besides, one has
 \beq\label{x0 nc}
x_0^{\text nc}(h) = \frac{27}{64}\, h^{-3} (1 + o(1)) \ \ \ \
\text{as} \ \ h \to +\infty.
 \eeq

The minimal resistance is calculated according to
(\ref{7}); after some algebra one gets
$$
\mathcal{R}_{\text nc}(h)\, =\, \frac{1}{2}\, +\, \frac{3 + 2z^{2/3}
- 8z^{4/3}}{16z^{5/3}}\, \sqrt{z^{2/3} - 1}\, +\, \frac{3}{16z^2}\,
\ln (z^{1/3} + \sqrt{z^{2/3} -1}).
$$
One also gets from theorem \ref{t2} that $\mathcal{R}_{\text
nc}(0) = 0.5$ and
 \beq\label{R nc}
 \mathcal{R}_{\text nc}(h) = \frac{27}{128} h^{-2} (1 +
o(1))  \ \ \ \ \text{as} \ \ h \to +\infty.
 \eeq
 \vspace{1mm}

{\bf 2.} The original {\it Newton problem} (case (i) in our
classification) is also a particular case of (\ref{1}), with $p(u) =
p_N(u) := 2/(1 + u^2)$. One has $u_0 = 1$ and $\bar p_N(u) = \left\{
\begin{array}{ll} ~ 2 - u & \text{if }\, 0 \le u \le 1\\ 2/(1 + u^2) &
\text{if }~~~ u \ge 1~~, \end{array} \right.$ and after some
calculation one gets that $\xi^N_0 = 1$ and the function
$\Upsilon_N(z)$,\, $z \ge 1$, in a parametric representation, is
$\Upsilon_N = \frac 14\, (3u^4/4 + u^2 - \ln u - 7/4)$,~ $z = (1 +
u^2)^2/(4u)$,\, $u \ge 1$. From here one obtains the well-known
Newton solution:~ if $0 \le x \le x_0$ then $f_h(x) = 0$, and if
$x_0 < x \le 1$ then $f_h$ is defined parametrically: $f_h =
\frac{x_0}{4}\, (3u^4/4 + u^2 - \ln u -  7/4)$,~ $x =
\frac{x_0}{4}\, \frac{(1 + u^2)^2}{u}$, where $x_0 = 4u_*/(1 +
u_*^2)^2$ and $u_*$ is determined from the equation $(3u_*^4/4 +
u_*^2 - \ln u_* - 7/4)\, u_*/(1 + u_*^2)^2 = h$. The function $f_h$
is not differentiable at $x_0$: one has $f'_h(x_0 + 0) = 1$ and
$f'_h(x_0 - 0) = 0$.

One also has $\mathcal{R}_N(0) = 1$,
 \beq\label{R nc}
 \mathcal{R}_{N}(h) = \frac{27}{32} h^{-2} (1 +
o(1))  \ \ \ \ \text{as} \ \ h \to +\infty.
 \eeq
 and
  \beq\label{x0 N}
x_0^{N}(h) = \frac{27}{16}\, h^{-3} (1 + o(1)) \ \ \ \ \text{as} \ \ h \to +\infty.
 \eeq
  \vspace{1mm}

{\bf 3.} The minimal problem in the {\it single impact case} with $h
> M^* \approx 0.54$ can also be reduced to (\ref{1}), with $p(u) =
p_{\text si}(u) := \left\{ \begin{array}{ll} ~ p^* & \text{if }\, u
= 0\\ 2/(1 + u^2) & \text{if }~~~ u > 0~~, \end{array} \right.$
where $p^* = 8(\ln(8/5) + \arctan(1/2) - \pi/4) \approx 1.186$. This fact
can be easily deduced from \cite{CL1}; for the reader's convenience
we put the details of derivation in the next section.\footnote{We
would like to stress that the results presented here about the
single impact case can be found in \cite{CL1} or can be easily
deduced from the main results of \cite{CL1}.} From the above formula
one can calculate that $u_0 \approx 1.808$ and $\xi_0^{\text si}
\approx 2.52$.

The asymptotic formulas here take the form
 \beq\label{x0 si}
x_0^{\text si}(h) = \xi_0^{\text si} \cdot x_{0}^N(h) (1 + o(1)) \ \
\ \ \text{as} \ \ h \to +\infty
 \eeq
and
 \beq\label{R si}
 \mathcal{R}_{\text si}(h) = \frac{27}{32} h^{-2} (1 +
o(1))  \ \ \ \ \text{as} \ \ h \to +\infty.
 \eeq
 Finally, using the results of
\cite{CL1}, one can show that $\mathcal{R}_{\text si}(0) = \pi/2 -
2\arctan(1/2) \approx 0.6435$. This will also be made in the next
section.
 \vspace{1mm}

Now we are in a position to compare the solutions in the three
cases. One obviously has $\mathcal{R}_{\text nc}(h) \le
\mathcal{R}_{\text si}(h) \le \mathcal{R}_N(h)$. From the above
formulas one sees that $\mathcal{R}_{\text nc}(0) = 0.5$,\,
$\mathcal{R}_N(0) = 1$, and $\mathcal{R}_{\text si}(0)  \approx
0.6435$. Besides, one has $\lim_{h \to +\infty} ({\mathcal{R}_{\text
nc}(h)}/{\mathcal{R}_{N}(h)}) = 1/4$ and $\lim_{h \to +\infty}
({\mathcal{R}_{\text si}(h)}/{\mathcal{R}_{N}(h)}) = 1$. Thus, for
"short" bodies, the minimal resistance in the nonconvex case is two
times smaller than in the Newton case, and $22\%$ smaller, as
compared to the single impact case. For "tall" bodies, the minimal
resistance in the nonconvex case is four times smaller as compared
th the Newton case, while the minimal resistance in the Newton case
and in the single impact case are (asymptotically) the same.

In the three cases of interest, the convex hull of the
three-dimensional optimal body of revolution has a flat disk of
radius $x_0(h)$ at the front part of its boundary. One always has
$x_0(0) = 1$. For "tall" bodies, one has $\lim_{h \to +\infty}
(x_0^{\text nc}(h) / x_{0}^N(h)) = 1/4$ and $\lim_{h \to +\infty}
(x_0^{\text si}(h) / x_{0}^N(h)) = \xi_0^{\text si} \approx 2.52$;
that is, the disk radius in the non-convex case and in the single
impact case is, respectively, 4 times smaller and $2.52$ times
larger, as compared to the Newton case.

Besides, in the nonconvex case, the front part of the surface of the
body's convex hull is smooth. On the contrary, in the Newton case,
the front part of the body's surface has singularity at the boundary
of the front disk.

\section{Proofs of the results}

\subsection{Proof of lemma \ref{l1}}

It suffices to show that
 \beq\label{l1-1}
 v^z_G(x) \ge \hat{v}^z_{\text{conv}\,G}(x) \ \ \text{for any} \ \ x \in [0,\, 1].
 \eeq

\begin{figure}
\begin{center}
\includegraphics[height=12cm, angle=-90]{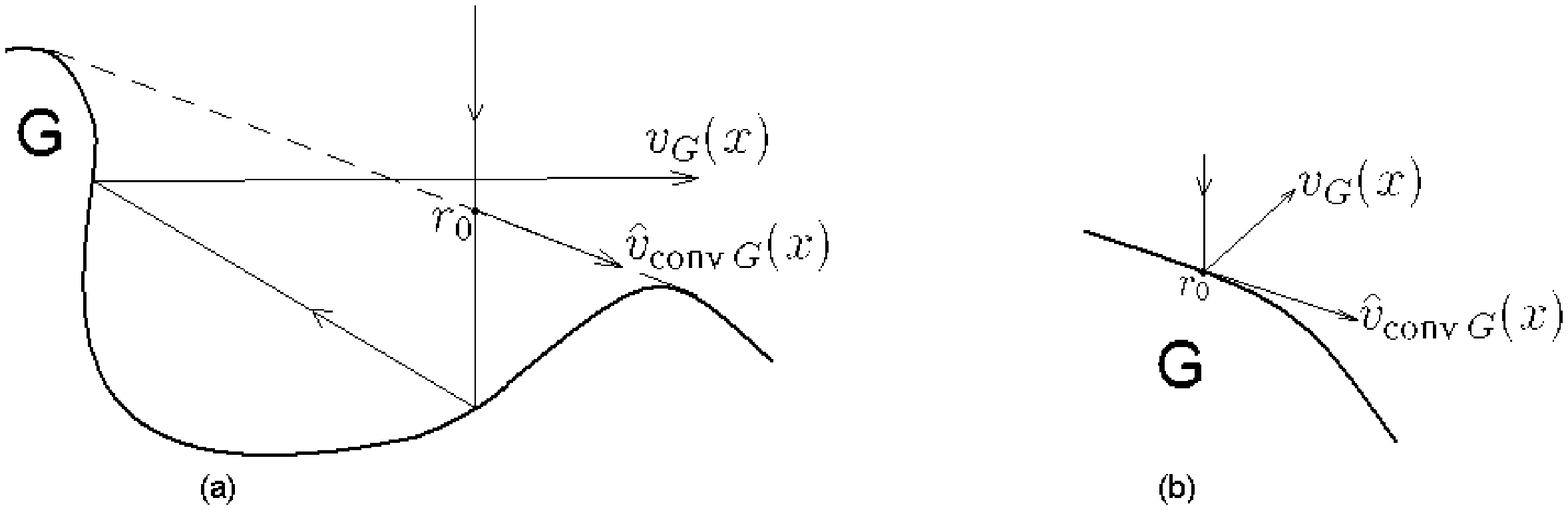}
\caption{Two scenarios of reflection.}
\label{conv}
\end{center}
\end{figure}

Consider two scenarios of motion for a particle that initially
moves vertically downwards, $x(t) = x$ and $z(t) = -t$. First, the
particle hits conv\,$G$ at a point $r_0 \in \pl(\text{conv}\, G)$
according to the modified reflection law and then moves with the
velocity $\hat{v}_{\text{conv}\,G}(x)$. Second, it hits $G$
(possibly several times) according to the law of elastic
reflection, and then moves with the velocity $v_G(x)$. Denote by
$n$ the outer unit normal to $\pl(\text{conv}\, G)$ at $r_0$; on fig.
\ref{conv} there are shown two possible cases: $\, r_0 \in \pl G$
and $r_0 \not\in \pl G$.

It is easy to see that
 \beq\label{l1-2}
 \langle v_G(x), n \rangle \ge 0,
 \eeq
where $\langle \cdot \,, \cdot \rangle$ means the scalar product.
Indeed, denote by $r(t) = (x(t), z(t))$ the particle position at
time $t$. At some instant $t_1$ the particle intersects
$\pl(\text{conv}\, G)$ and then moves outside $\text{conv}\, G$. The
function $\langle r(t), n \rangle$ is linear and satisfies $\langle
r(t), n \rangle \ge \langle r(t_1), n \rangle$ for $t \ge t_1$,
therefore its derivative $\langle v_G(x), n \rangle$ is positive.

>From (\ref{l1-2}) and the relations $\langle \hat
v_{\text{conv}\,G}(x), n \rangle = 0$,\, $\hat
v^z_{\text{conv}\,G}(x) \le 0$ and $n_z \ge 0$,\, $n_x \ge 0$ one gets (\ref{l1-1}).

\subsection{Proof of lemma \ref{l2}}

Take a family of piecewise affine even functions $f_\ve: [-1,\, 1]
\to [0,\, h]$ such that $f_\ve'$ uniformly converges to $f_G'$ as
$\ve \to 0^+$. Require also that the functions $f_\ve$ are concave
and monotone decreasing as $x > 0$, and $f_\ve(0) = h$,\, $f_\ve(1)
= f_G(1)$. Consider the family of convex sets $G_\ve \in
\GGG_h^{\text conv}$ bounded from above by the graph of $f_\ve$ and
from below, by the segment $-1 \le x \le 1$,\, $z = 0$. Taking into
account (\ref{resmod}), one gets $\lim_{\ve \to 0^+} \widehat
R(G_\ve) = \widehat R(G)$.

Below we shall determine a family of sets $G_{\ve,\del} \in \GGG_h$
such that $\lim_{\del \to 0^+} R(G_{\ve,\del}) = \widehat R(G_\ve)$
and next, using the diagonal method, select a sequence
$\ve_n \to 0, \ \del_n \to 0$ such that $\lim_{n\to\infty} R(G_{\ve_n,\del_n})
= \lim_{n \to \infty} \widehat R(G_{\ve_n}) = \widehat R(G)$. This will
finish the proof.

Fix $\ve > 0$ and denote by $-1 = x_{-m} < x_{-m+1} < \ldots < x_0 =
0 < \ldots < x_m = 1$ the jump values of the piecewise constant
function $f'_\ve$. (One obviously has $x_{-i} = -x_i$.) For each $i
= 1, \ldots, m$ we shall define a non self-intersecting curve
$l^{i,\ve,\del}$ that connects the points $(x_{i-1},
f_\ve(x_{i-1}))$ and $(x_{i}, f_\ve(x_{i}))$ and is contained in the
quadrangle $x_{i-1} \le x \le x_{i}$,\, $f_\ve(x_i) \le z \le
f_\ve(x_{i-1}) + (f'_\ve(x_{i-1}+0)+\delta) \cdot (x - x_{i-1})$.
The curve $l^{-i,\ve,\del}$ is by definition symmetric to
$l^{i,\ve,\del}$ with respect to the axis $Oz$. Let now
$l^{\ve,\del} := \cup_{-m \le i \le m} l^{i,\ve,\del}$ and let
$G_{\ve,\del}$ be the set bounded by the curve $l^{\ve,\del}$, by
the two vertical segments $0 \le z \le f_\ve(1)$,\, $x = \pm 1$, and
by the horizontal segment $-1 \le x \le 1$,\, $z = 0$.

For an interval $I \subset [0,\, 1]$, define
 \beq\lab{rihat}
\widehat R_I(G_\ve) := \int_I (1 + \hat v_{G_\ve}^z(x))\, x\, dx
 \eeq
and
 \beq\lab{ri}
 R_I(G_{\ve,\del}) := \int_I (1 + v_{G_{\ve,\del}}^z(x))\, x\, dx.
 \eeq
Denote $I_i = [x_{i-1},\, x_i]$; one obviously has $\widehat
R(G_\ve) = \sum_{i=1}^m \widehat R_{I_i}(G_\ve)$ and
$R(G_{\ve,\del}) = \sum_{i=1}^m R_{I_i}(G_{\ve,\del})$. Thus, it
remains to determine the curve $l^{i,\ve,\del}$ and prove that
\begin{equation}\label{l2 convergence}
\lim_{\del \to 0^+} R_{I_i}(G_{\ve,\del}) = \widehat R_{I_i}(G_\ve).
\end{equation}
This will complete the proof of the lemma.

Note that for $x \in I_i$,\, $i = 1, \ldots, m$ holds
 \beq\label{v hat}
\hat v^z_{G_{\ve}} = \frac{f'_\ve(x_{i-1} + 0)}{\sqrt{1 +
(f'_\ve(x_{i-1} + 0))^2}}.
 \eeq

\begin{figure}
\begin{center}
\includegraphics[width=8cm, angle=-90]{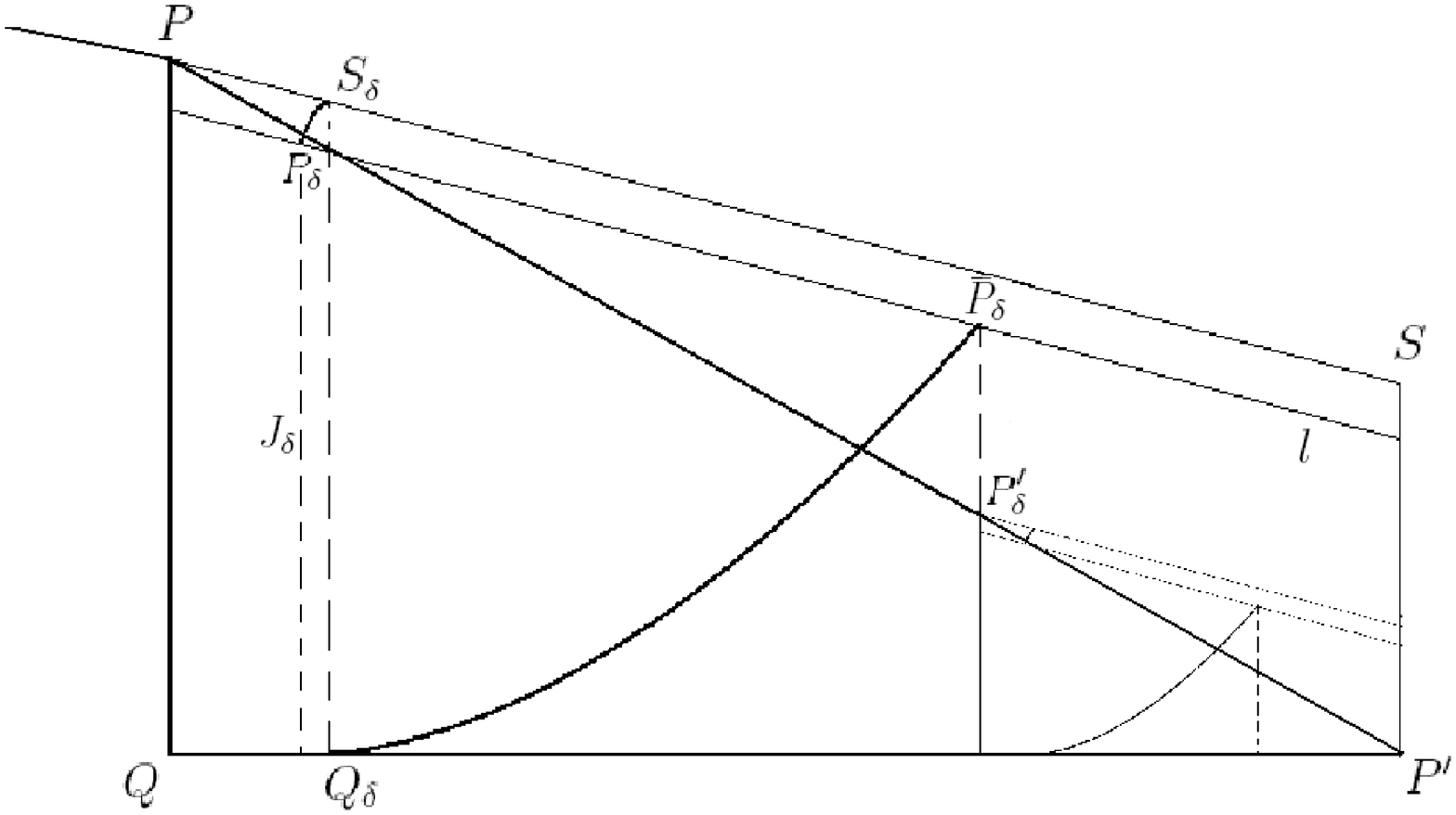}
\caption{Constructing the curve $l^{i,\ve,\del}$: a detailed view.}
\label{constr1}
\end{center}
\end{figure}

Fix $\ve$ and $i$ and mark the points $P = (x_{i-1},
f_\ve(x_{i-1}))$,\, $P' = (x_{i}, f_\ve(x_{i}))$,\, $Q = (x_{i-1},
f_\ve(x_{i}))$, and $S = (x_i,\, f_\ve(x_{i-1}) +
(f'_\ve(x_{i-1}+0)+\delta) \cdot (x_i - x_{i-1}))$; see
fig.\,\ref{constr1}. Mark also the point $Q_\del = (x_{i-1} + \del,
f_\ve(x_{i}))$, which is located on the segment $QP'$ at the
distance $\del$ from $Q$, and the points $P_\del = (x_{i-1} + \del,
f_\ve(x_{i-1} + \del))$ and $S_\del = (x_{i-1} + \del,
f_\ve(x_{i-1}) + (f'_\ve(x_{i-1}+0)+\delta) \cdot \del)$, which have
the same abscissa as $Q_\del$ and belong to the segments $PP'$ and
$PS$, respectively. Denote by $l$ the line that contains $P_\del$
and is parallel to $PS$. Denote by $\Pi_\del$ the arc of the
parabola with vertex $Q_\del$ and focus at $P_\del$ (therefore its
axis is the vertical line $Q_\del P_\del$). This arc is bounded by
the point $Q_\del$ from the left, and by the point $\bar P_\del$ of
intersection of the parabola with $l$, from the right. Denote by
$x_i^\del$ the abscissa of $\bar P_\del$ and denote by $P_\del'$ the
point that lies in the line $PP'$ and has the same abscissa
$x_i^\del$. Denote by $\pi_\del$ the arc of the parabola with the
same focus $P_\del$, the axis $l$, and the vertex situated on $l$ to
the left from $P_\del$. The arc $\pi_\del$ is bounded by the vertex
from the left, and by the point $S_\del'$ of intersection of the
parabola with the line $Q_\del P_\del$, from the right. There is an
arbitrariness in the choice of the parabola; let us choose it in
such a way that the arc $\pi_\del$ is situated below the line $PS$.
Finally, denote by $J_\del$ the perpendicular dropped from the left
endpoint of $\pi_\del$ to $QP'$, and denote by $Q_\del'$ the base of
this perpendicular.

If $x_i^\del \ge x_i$, the curve $l^{i,\ve,\del}$ is the union
(listed in the consecutive order) of the segments $PS_\del$ and
$S_\del S_\del'$, the arc $\pi_\del$, the segments $J_\del$ and
$Q_\del' Q_\del$, and the part of $\Pi_\del$ located to the left of
the line $P'S$.

If  $x_i^\del < x_i$, the definition of $l^{i,\ve,\del}$ is more
complicated. Define the homothety with the center at $P'$ that sends
$P$ to $P_\del'$, and define the curve $\tilde l^{i,\ve,\del}$ by
the following conditions:~ (i) the intersection of $\tilde
l^{i,\ve,\del}$ with the strip region $x_{i-1} \le x \le x_{i}^\del$
is the union of $PS_\del$,\, $S_\del S_\del'$,\, $\pi_\del$,\,
$J_\del$,\, $Q_\del' Q_\del$,\, $\Pi_\del$, and the interval $\bar
P_\del P_\del'$;~ (ii) under the homothety, the curve $\tilde
l^{i,\ve,\del}$ moves into itself. The curve $\tilde l^{i,\ve,\del}$
is uniquely defined by these conditions; it does not have
self-intersections and connects the points $P$ and $P'$. However, it
is not piecewise smooth, since it has infinitely many singular
points near $P'$. In order to improve the situation, define the
piecewise smooth curve $l^{i,\ve,\del}$ in the following way: in the
strip $x_{i-1} \le x < x_i - \del$, it coincides with $\tilde
l^{i,\ve,\del}$, the intersection of $l^{i,\ve,\del}$ with the strip
$x_i - \del < x \le x_i$ is the horizontal interval $x_i - \del < x
\le x_i$,\ $z = f_\ve(x_i)$, and the intersection of
$l^{i,\ve,\del}$ with the vertical line  $x = x_i - \del$ is a point
or a segment (or maybe the union of a point and a segment) chosen in
such a way that the resulting curve $l^{i,\ve,\del}$ is continuous.

\begin{figure}\label{lastpic}
\begin{center}
\includegraphics[width=8cm, angle=-90]{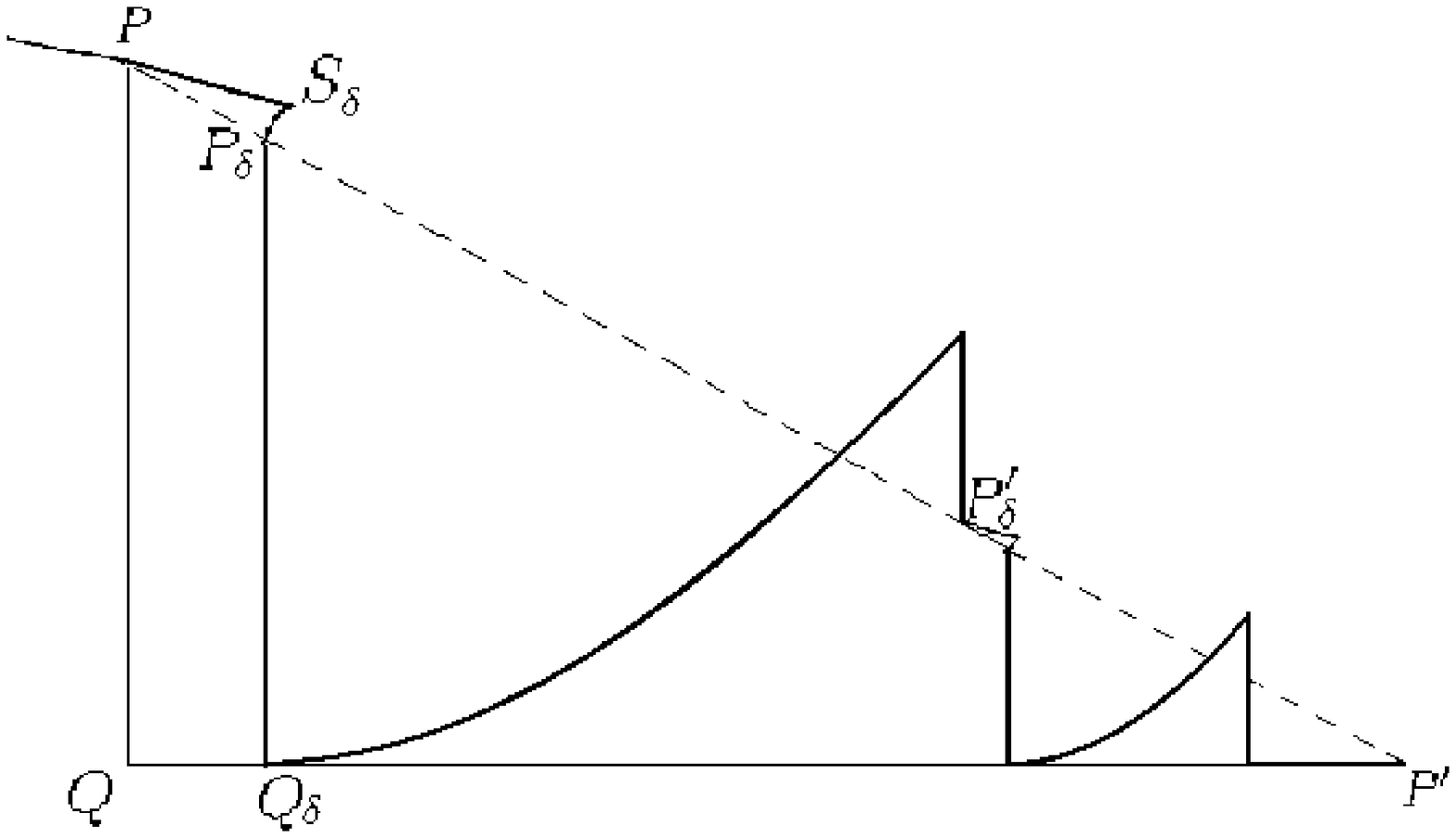}
\caption{The curve $l^{i,\ve,\del}$, again.}
\label{constr2}
\end{center}
\end{figure}

The particles of the flow falling on the arc $\Pi_\del$ make a
reflection from it, pass through the focus $P_\del$, then make
another reflection from the arc $\pi_\del$, and finally move freely, the velocity being
parallel to $l$. Choose $\del < |f_\ve'(0^+)|$ and $\del < \min_{1
\le i \le m-1} (f_\ve'(x_{i-1}+0) - f_\ve'(x_{i}+0))$, then the
particles after the second reflection will never intersect the other
curves $l^{j,\ve,\del}$,\, $j \ne i$. Thus, for the corresponding
values of $x$, the vertical component of the velocity of the
reflected particle is
 \beq\label{l2 reflected}
v^z_{G_{\ve,\del}}(x) = \frac{f'_\ve(x_{i-1} + 0) + \del}{\sqrt{1 +
(f'_\ve(x_{i-1} + 0) + \del)^2}} = \hat v^z_{G_\ve}(x) + O(\del),
~~~ \del \to 0^+.
 \eeq
If $x_i^\del \ge x_i$, the formula (\ref{l2 reflected}) is valid for
$x \in [x_{i-1} + \del,\, x_i]$. If $x_i^\del < x_i$, it is valid
for the values $x \in [x_{i-1} + \del,\, x_i^\del]$. Note, however,
that (\ref{l2 reflected}) is also valid for values of $x$ that
belong to the iterated images of $x \in [x_{i-1} + \del,\,
x_i^\del]$ under the homothety, but do not belong to $[x_i - \del,\,
x_i]$. Summarizing, (\ref{l2 reflected}) is true for $x \in
[x_{i-1},\, x_i]$, except for a set of
values of measure $O(\del)$. Thus, taking into account (\ref{rihat}),\, (\ref{ri}),\,
(\ref{v hat}), and (\ref{l2 reflected}), the convergence (\ref{l2
convergence}) is proved. Q.E.D.

\subsection{Proof of theorem \ref{t2}}

Let us first state (without proof) the following lemma.

\begin{lemma}\label{pl}
Let $\lam > 0$ and let the function $f \in {\mathcal F}_h$ satisfy
the condition

{$\mathbf I_\lam$}. \quad $f(1) = h$, and for almost all $x \in
[0,\, 1]$ the value $u = f'(x)$ is a solution of the problem
\begin{equation}\label{a}
x p(u)+\lam u \rightarrow  \min, \quad u \in \RRR_+.
\end{equation}
Then the function  $f$ is a solution of the problem (\ref{1}) and
any other solution satisfies the condition  $\mathrm{I}_\lam$ with
the same value of $\lam$.
\end{lemma}

This simple lemma is a direct consequence of the Pontryagin maximum
principle. Its proof can be found, for example, in \cite{T} or in
\cite{temperature}.

Now we shall find the function $f_h$ satisfying the condition
$\mathrm I_\lam$ for some positive $\lam$. Let  $x \in [0,\, 1]$ be
the value for which $\mathrm I_\lam$ is fulfilled. Then the value $u
= f'_h(x)$ is also a minimizer for the function $x \bar p(u)+\lam
u$, and $p(u) = \bar p(u)$. This implies that (if the function $f_h$
really exists then)
 \beq\label{rfh}
 \mathcal{R}[f_h] = \int_0^1 \bar p(f'_h(x))\, x\, dx.
 \eeq
Besides, if $u > 0$ and $\bar p$ is differentiable at $u$ then one
has $\frac{d}{du} (x \bar p(u) + \lam u) = 0$, hence
 \beq\label{3}
\frac{x}{\lam}\ =\ -\frac{1}{\bar p'(u)}\,.
 \eeq
If $u > 0$ and $\bar p$ is not differentiable at $u$, then it has
left and right derivatives at this point and
 \beq\label{3.1}
-\frac{1}{\bar p'(u-0)}\ \le\ \frac{x}{\lam}\ \le\ -\frac{1}{\bar
p'(u+0)}\,.
 \eeq
If, finally, $u = 0$ then one has
  \beq\label{3.2}
 \frac{x}{\lam}\ \le\ -\frac{1}{\bar
p'(0)}\ =\ \xi_0.
 \eeq
Put $z = 1/\lam$ and $x_0 = \xi_0/z$ and rewrite (\ref{3}) and
(\ref{3.1}) in terms of the generalized inverse function: $\upsilon(zx-0) \le u \le \upsilon(zx)$; thus the equality
 \beq\label{4}
u\ =\ \upsilon(z x),
 \eeq
is valid for almost all values $x \ge x_0$. Taking
into account (\ref{3.2}), substituting $u = f'_h(x)$, and
integrating both parts of (\ref{4}) with respect to $x$, one comes
to (\ref{f solution}). In particular, $f_h'(x_0 + 0) =
\upsilon(\xi_0 + 0) = u_0$. Using that $f_h(1) = h$, one gets
(\ref{z solution}).

The function $\Upsilon(z)/z$ is continuous and monotone increasing;
it is defined on $[\xi_0, +\infty)$ and takes the values from 0 to
$+\infty$. Therefore the equation (\ref{z solution}) uniquely
defines $z$ as a continuous monotone increasing function of $h$; in
particular, $z(0) = \xi_0$ and $x_0(0) = \xi_0/z(0) = 1$. The
relations (\ref{f solution}) and (\ref{z solution}) define the
function $f_h$ solving the minimization problem (\ref{1}). From the
construction one can see that this function is uniquely defined.

Recall that $\mathcal{R}(h) = \mathcal{R}[f_h]$. Integrating by
parts the right hand side of (\ref{rfh}), one gets
$$
\mathcal{R}(h)\ =\ \frac{\bar p(f_h'(1))}{2} \ -\ \int_0^1\, \frac{\
x^2}{2}\ \bar p'(f_h'(x))\, df_h'(x).
$$
Taking into account that  $f_h'(1) = \upsilon(z)$ and $x \bar
p'(f_h'(x)) = -\lam = -1/z$, one obtains
$$
\mathcal{R}(h)\ =\ \frac{\bar p(\upsilon(z))}{2}\, +\,
\frac{1}{2z}\, \int_0^1 x\, df_h'(x),
$$
and integrating by parts once again, one gets (\ref{7}).
Substituting in (\ref{7}) $h = 0$ and using that $z(0) = \xi_0$,\,
$\upsilon(\xi_0) = u_0$,\, $\Upsilon(\xi_0) = 0$, one obtains
$\mathcal{R}(0) = (\bar p(u_0) + u_0/\xi_0)/2$, and using that $p(0)
- \xi_0^{-1} u_0 = \bar p(u_0)$, one obtains $\mathcal{R}(0) =
p(0)/2$.

Taking into account the asymptotic of $\bar p\,$ (which is the same
as the asymptotics of $p$:~ $\bar p(u) = c\, u^{-\al} (1 + o(1))$,\,
$u \to +\infty$), and the asymptotic of $\bar p'$:~ $\bar p'(u) =
-c\al\, u^{-\al - 1} (1 + o(1))$,\, $u \to +\infty,\,$ one comes to
the formulas
$$
\upsilon(\xi)\, =\, (C\al)^{\frac{1}{\al + 1}}\, \xi^{\frac{1}{\al +
1}}\, (1 + o(1)),~~~ \xi \to +\infty,
$$
$$
\Upsilon(z)\, =\, \left( \frac{\al + 1}{\al + 2}
\right)(c\al)^{\frac{1}{\al + 1}}\, z^{\frac{\al + 2}{\al + 1}}\, (1
+ o(1)),~~~ z \to +\infty,
$$
and
$$
z\, =\, \frac{1}{c\al} \left( \frac{\al + 2}{\al + 1} \right)^{\al +
1}\, h^{\al + 1}\, (1 + o(1)),~~~ h \to +\infty.
$$
Substituting them into (\ref{7}) and using the relation $x_0 =
\xi_0/z$, after a simple algebra one obtains (\ref{8}) and
(\ref{9}). The theorem is proved.
 \vspace{2mm}

Summarizing, the three-dimensional bodies of revolution minimizing
the resistance are constructed as follows. First, we find the
function $f_h^{\text nc}$ minimizing the functional (\ref{probl})
and define the convex set $-1 \le x \le 1$,\, $0 \le z \le h -
f_h^{\text nc}(|x|)$. Next, the upper part of its boundary (which is
the graph of the function $z = h - f_h^{\text nc}(|x|)$) is
approximated by a broken line and then substituted with a curve with
rather complicated behavior, according to lemma \ref{l2}. The set
bounded from above by this curve is "almost convex": it can be
obtained from a convex set by making small hollows on its boundary.
By rotating it around the axis $Oz$, one obtains the body of
revolution $B$ having nearly minimal resistance $\mathrm{R}_z(B)$.

The vertical central cross sections of optimal bodies in the Newton,
single impact, and nonconvex cases, for $h = 0.8$, are presented
on figure 9. 

\begin{figure}\label{scat3}
\begin{center}
\includegraphics[width=8cm, angle=0]{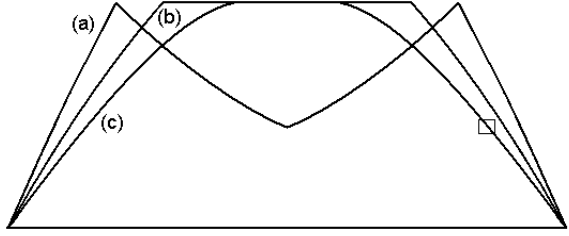}
\caption{Profiles of optimal solutions in the single impact (a),
Newton (b), and nonconvex (c) cases, for $h=0.8$. In the nonconvex
case, the profile is actually a zigzag curve with very small
zigzags, as shown on the next figure. } \label{normRF}
\end{center}
\end{figure}

\begin{figure}
\begin{center}
\includegraphics[height=7cm, angle=-90]{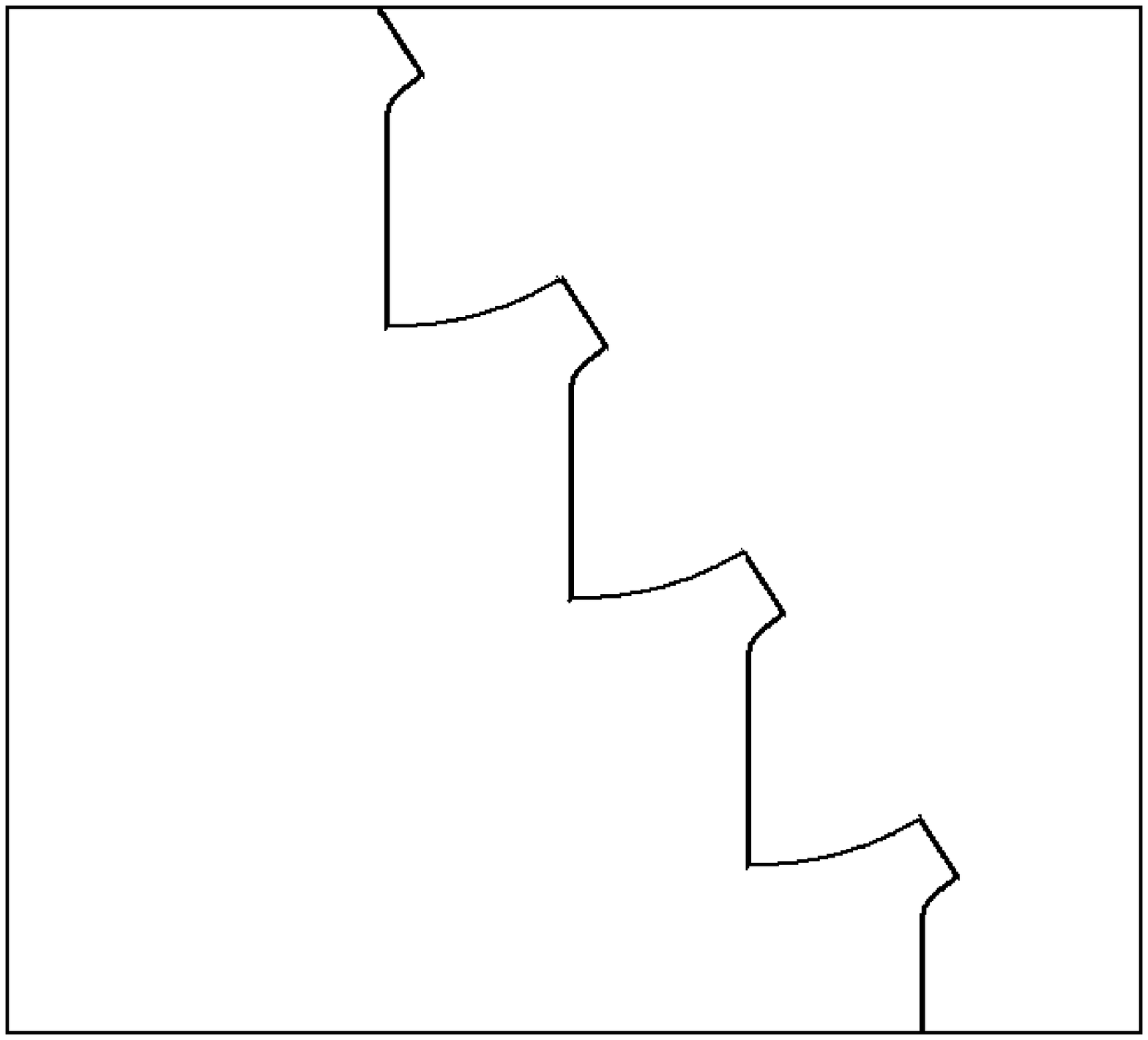}
\caption{Detailed view of the zigzag curve.}
\label{micr}
\end{center}
\end{figure}

\subsection{Derivation of the asymptotic relations in the single impact case}

For $h$ small (namely, $h < M^* \approx 0.54$), a solution in the
single impact case can be described as follows. There are marked
several values $-1 < x_{-2n+1} < x_{-2n+2} < \ldots < x_{2n-2} <
x_{2n-1} < 1$,\, $n \ge 2$ related to the singular points of the
solution. As $h \to 0^+,\ n = n(h)$ goes to infinity. One has
$x_{-k} = -x_k$ and $x_{2i} = (x_{2i-1} + x_{2i+1})/2$; thus $x_0 =
0$. Besides, one has $\max_k(x_k - x_{k-1}) = x_1 = 4h/3$. The
vertical central cross section of the solution $G = G^{\text si}_h
\subset \RRR^2_{x,z}$ is bounded from above by the graph of a
continuous non-negative piecewise smooth even function $f =
f_h^{\text si}$, and from below, by the segment $-1 \le x \le 1$,\,
$z = 0$. This function has singularities at the points $x_k$, and
the values of the function at the points $x_{2i-1}$ coincide:
$f(x_{2i-1}) = h$. On each interval $[x_{2i-1},\, x_{2i}]$, the
graph of $f$ is the arc of parabola with vertical axis and with the
focus at $(x_{2i+1}, h)$. Similarly, on $[x_{2i},\, x_{2i+1}]$ the
graph of $f$ is the arc of parabola with vertical axis and with the
focus at $(x_{2i-1}, h)$. The first parabola contains the focus of
the second one, and vice versa. From this description one can see
that on $[x_{2i-1},\, x_{2i}]$, the function equals $f(x) = \frac{(x
- x_{2i+1})^2}{2(x_{2i+1} - x_{2i-1})} + y_i$, and on $[x_{2i},\,
x_{2i+1}]$,~ $f(x) = \frac{(x - x_{2i-1})^2}{2(x_{2i+1} - x_{2i-1})}
+ y_i$, where $y_i = h - (x_{2i+1} - x_{2i-1})/2$. On the intervals
$[-1,\, x_{-2n+1}]$ and $[x_{2n-1},\, 1]$ the graph of the function
represents the so-called "Euler part" of the solution (see
\cite{CL1}).

Note that the solution is not unique. The values $x_1$ and
$x_{2n-1}$ are uniquely determined, but there is arbitrariness in
choice of the intermediate values $x_3, \ldots, x_{2n-3}$ and also
in the number $n$ of the independent parameters.

After some calculation, one obtains the value of $1 + v_G^z$ for the
figure $G$:
$$
\text{if}~~ x \in [x_{2i-1},\, x_{2i}], \ \ \ \ \ \ \ 1 + v_G^z(x)\
=\ \frac{2}{1 + \left( \frac{x_{2i+1} - x}{x_{2i+1} - x_{2i-1}}
\right)^2}~;
$$
$$
\text{if}~~ x \in [x_{2i},\, x_{2i+1}], \ \ \ \ \ \ \ 1 + v_G^z(x)\
=\ \frac{2}{1 + \left( \frac{x - x_{2i-1}}{x_{2i+1} - x_{2i-1}}
\right)^2}~.
$$

Let us now calculate the integral $\int_{x_{2i-1}}^{x_{2i+1}} (1 +
v_G^z(x))\, x\, dx$,\, $1 \le i \le n-1$. Since the function $1 +
v_G^z(x)$,\, $ x \in [x_{2i-1},\, x_{2i+1}]$ is symmetric with
respect to $x = x_{2i}$, the integral equals $2x_{2i}
\int_{x_{2i}}^{x_{2i+1}} (1 + v_G^z(x))\, dx$. Changing the variable
$t = (x - x_{2i-1})/(x_{2i+1} - x_{2i-1})$ and taking into account
that $2x_{2i} = x_{2i-1} + x_{2i+1}$, one comes to the integral
$2x_{2i} (x_{2i+1} - x_{2i-1}) \int_{1/2}^1 2/(1 + t^2)\, dt =
(x_{2i+1}^2 - x_{2i-1}^2) (\pi/2 - 2\arctan(1/2))$. Therefore
$$
\int_{x_{1}}^{x_{2n-1}} (1 + v_G^z(x))\, x\, dx\ =\ (x_{2n-1}^2 -
x_{1}^2) (\pi/2 - 2\arctan(1/2)).
$$
Taking into account that $x_1 = 4h/3 \to 0$ and $x_{2n(h)-1} \to 1$
as $h \to 0^+$, one finally gets
$$
R(G^{\text si}_h)\, =\, \int_0^1 (1 + v_{G^{\text si}_h}^z(x))\, x\,
dx\, =\, \pi/2 - 2\arctan(1/2) + o(1), \ \ \ \ \ h \to 0^+,
$$
that is, $\mathcal{R}_{\text si}(0) = \pi/2 - 2\arctan(1/2) \approx
0.6435$.

If $h > M^*$, the function $f = f_h^{\text si}$ has three singular
points: $x_1 = x_1(h),\ 0$, and $-x_1$. On the interval $[-x_1,\,
x_1]$, the graph of $f$ is the union of two parabolic arcs, as
described above with $i = 0$. On the intervals $[-1,\, -x_1]$ and
$[x_1,\, 1]$, the graph is the "Euler part" of the solution; on both
intervals, $f$ is a concave monotone function, with $f(\pm 1) = 0$
and $f(\pm x_1) = h$. The part of resistance of $G = G^{\text si}_h$
related to $[0,\, x_1]$ can be calculated:
$$
\int_{0}^{x_{1}} (1 + v_G^z(x))\, x\, dx\ =\ x_1 p^*,
$$
where $p^* = 8(\ln(8/5) + \arctan(1/2) - \pi/4) \approx 1.186$. That
is, the convex hull of $G$ represents the solution of the problem
(\ref{1}) with $p(u) = p_{\text si}(u) = \left\{ \begin{array}{ll} ~
p^* & \text{if }~ u = 0\\ 2/(1 + u^2) & \text{if }~ u > 0~~.
\end{array} \right.$


\end{document}